# REMARKS ON METHODS OF FONTAINE AND FALTINGS

KIRTI JOSHI

1. INTRODUCTION

Let $G = \text{Gal}(\bar{\mathbf{Q}}/\mathbf{Q})$ be the absolute Galois group. Let $p$ be a prime and let $\rho : G \to GL_2(\mathbf{F}_p)$ be an absolutely irreducible Galois representation. Under certain conditions, viz., that $\rho$ is odd, by [14], one expects $\rho$ to arise from a modular form of a suitable level, weight and nebentype. Serre also gives a recipe for computing these optimal invariants attached to $\rho$, which will be denoted by $N(\rho)$ (level), $k(\rho)$ (weight) and $\epsilon(\rho)$ (nebentype).

The situation of interest to us is the following. Let us assume that $\rho$ is odd, unramified outside $p$ with Serre weight $k(\rho) = 2$ (so $\rho$ is finite flat at $p$). By Serre's conjecture one expects that $\rho$ is modular of some level, weight and nebentype. Then by level lowering results of Mazur, Ribet and others, (see [13]) such a $\rho$ must arise from a modular form of weight two and level one. But there are no such forms. Hence there can be no such $\rho$.

Thus one might ask if one can prove such non-existence results unconditionally (i.e. without Serre's conjecture). In this note we show that there no irreducible $\rho$ of the type considered in the above paragraph for small primes $p$. Our method, which is based on Fontaine's work also yields non-existence statements for higher weight representations. We also present a small generalization of Fontaine's method which allows auxiliary, but controlled, ramification.

Recall that there is a natural two dimensional irreducible Galois representation of Hodge-Tate weight $(0, 11)$ arising from the Ramanajun modular form $\Delta$ of level one. This representation is crystalline at 11 and unramified outside 11. So weight eleven is a natural boundary for these methods. However one could make the Fontaine-Mazur finiteness conjectures a little more explicit in this case and conjecture that for all but finitely many primes $p$, there is up to isomorphism, only one irreducible two dimensional $\mathbf{Q}_p$-representation of $G$ which of Hodge-Tate weight $(0, 11)$ and which is not a Tate twist of a lower weight Hodge-Tate representation and which is crystalline at $p$ and unramified else where. In fact we may even ask this question for the corresponding mod $p$ representation provided that we avoid the finite of of primes $p$ for which the mod $p$ Ramanujan representation is reducible.

In this rubric, as was pointed out to us by C. Khare, one notes that Serre's conjecture implies certain boundedness results for mod $p$ Galois representations. To see this we consider odd, irreducible representations $\rho : \text{Gal}(\bar{\mathbf{Q}}/\mathbf{Q}) \to GL_2(\bar{\mathbf{F}}_p)$ with bounded Artin conductor. Serre's conjecture implies that these representations arise from modular forms of a bounded level and weight and the space of such forms has bounded dimension and hence the set of isomorphism classes of such representations is finite (in fact we can make this more precise: a mod $p$ modular representation of level $N$ also arises from $S_2(N'p^2)$ where $N'$ is prime to $p$; so if the Artin conductor is bounded and $p$ is fixed we deduce boundedness). This is a modular variant of the Fontaine-Mazur finiteness conjectures (see [10]).

In the last section of the paper we note that an elementary variation of Faltings' method of proof of Shafarevitch conjecture yields a finiteness results for algebro-geometric representations (see Theorem 5.1). To be more specific we note that the set of isomorphism classes of irreducible representations of fixed dimension of the Galois group of a number field (or a function field), which are of fixed weight (in the sense of Deligne) and which are unramified outside a fixed finite set of primes, is finite. This result, whose proof is adapted almost verbatim from Faltings' proof, can viewed as the archimedean analogue of the Fontaine-Mazur conjecture. In the Fontaine-Mazur conjecture, stronger assumptions are made about the representations at primes dividing $p$ and conjecturally this information appears to determine the archimedean properties of the eigenvalues of Frobenius at almost all primes. As a consequence of the result mentioned above we deduce that there are only finitely many irreducible two dimensional $p$-adic representations which look like the representation arising from the Ramanujan motive of weight 11.

After this note was written, Professor Nigel Boston pointed out the work of Sharon Brueggeman, John Jones and David Roberts and we are indebted to him for these references. In her Ph.D thesis, Sharon Brueggeman (see [1] has obtained strong results on the image of Galois representations with values in $GL_2(\bar{\mathbf{F}}_5)$, which are unramified out side 5, under the assumption of Generalized Riemann Hypothesis. In [11], John Jones and David Roberts classify and give a complete list of sextic number fields whose discriminants are divisible only by the primes in $\{2,3\}$.

I am deeply indebted to Minhyong Kim for sharing with me his insights and for many discussions on matters pertinent to this work. We would like to thank C. Khare, Yihsiang Liow, Klaus Lux, Arvind Nair, C. .S Rajan, Dinesh Thakur and D. Ulmer for discussions.

## 2. Representations of Hodge-Tate Weight one

Let $\rho$ be a representation of $G = \text{Gal}(\bar{\mathbf{Q}}/\mathbf{Q})$, with values in $GL_m(\mathbf{F}_p)$. We will use the terminology of Fontaine (see [6]). We will say $\rho$ is crystalline at $p$ of Hodge-Tate weight at most $p-1$ if the restriction of $\rho$ to the decomposition group at $p$ is crystalline of Hodge-Tate weights at most $p-1$ in the sense of [9]. In other words the restriction of $\rho$ to the decomposition group at $p$ arises from a Fontaine-Laffaille module of filtration length at most $p-1$. The length of the filtration of the associated Fontaine-Laffaille module will be called the Hodge-Tate weight of the representation $\rho$. In particular $\rho$ is crystalline at $p$ of weight one if and only if $\rho$ is finite flat at $p$ and has Serre weight $k(\rho) = 2$.

We need the following variant of a theorem of Fontaine [7]. This theorem gives a very strong upper bound on $|d_K|^{1/n}$. On the other hand by methods of Odlyzko (see [12]), we get lower bounds for this quantity. These two bounds together give an upper bound on $n$.

**Theorem 2.1.** *Let $m \geq 2$ be an integer, let $S$ be a finite set of primes not containing $p$. Let $\rho : G \to GL_m(\mathbf{F}_p)$ be a representation with following properties:*
  1. *$\rho$ is unramified outside $S \cup \{p\}$*
  2. *$\rho$ is crystalline of Hodge-Tate weight $r \leq p-1$ at $p$.*
  3. *and $\rho$ is semi-stable at primes in $S$.*

Let $K$ be the fixed field of $\ker(\rho)$ and $n = [K : \mathbf{Q}]$ and let $d_K$ be discriminant of $K/\mathbf{Q}$. Then

$$(2.1) \qquad |d_K|^{1/n} < \left(\prod_{q \in S} q\right) p^{1 + \frac{r}{p-1}}.$$

*Proof.* The contribution of the $p$ part comes from Theorem 2 of [8]. For the contribution from a prime $q \in S$ we note that our assumption that $\rho$ is semi-stable at $q \in S$ implies that the image of the inertia subgroup of any prime lying over $q$ is a unipotent subgroup of $GL_m(\mathbf{F}_p)$, and hence has order coprime to $q$. Thus the ramification at any prime lying over $q$ in $K$ is tame. So now we use the standard bound on the different of a tamely ramified extension provided by the Lemma below.

The rest of the proof follows from the proof of Theorem 3 of [7]. □

**Lemma 2.2.** *Let $F/\mathbf{Q}_q$ be any tamely ramified finite extension of a $q$-adic field. Normalize the valuation $v$ on $F$ by setting $v(q) = 1$, and let $e$ be the ramification degree of $F/\mathbf{Q}_q$, $\delta_{F/\mathbf{Q}}$ be the different of $F$ and $d_{F/\mathbf{Q}_q}$ the discriminant of $F/\mathbf{Q}_q$. Then*

$$v(\delta_{F/\mathbf{Q}}) < 1 - \frac{1}{e} \leq 1 - \frac{1}{[F : \mathbf{Q}_q]} < 1.$$

*In particular if $K/\mathbf{Q}$ is a finite Galois extension which is tamely ramified at a prime $q$ and is unramified elsewhere then*

$$|d_{K/\mathbf{Q}}|^{1/[K:\mathbf{Q}]} < q.$$

*Proof.* This follows from Proposition 13 (page 58) of [15], note that Serre's normalization of $v$ is different from ours, which is consistent with Fontaine's normalization. □

**Remark 2.3.** Let $\rho : \mathrm{Gal}(\bar{\mathbf{Q}}/\mathbf{Q}) \to GL_2(\mathbf{F}_p)$ be any representation. Then any prime $q > p+1$ does not divide the order of $GL_2(\mathbf{F}_p)$ and hence $\rho$ is tamely ramified at such a prime. In particular the contribution of such a prime to discriminant of $K$ is given by the above lemma. Thus for any such $\rho$ we need only understand the contribution from primes $q$ such that $2 \leq q \leq p+1$.

We now assume that $\rho : \mathrm{Gal}(\bar{\mathbf{Q}}/\mathbf{Q}) \to GL_m(\mathbf{F}_p)$ is unramified outside $p$ and $\rho$ is finite flat at $p$. Let $K$ be the fixed field of $\ker(\rho) \subset G$. Then $K/\mathbf{Q}$ is a finite extension which is unramified outside $p$. As $\rho$ is finite flat at $p$ and unramified outside $p$ we see that $\det(\rho) = \chi$ the standard cyclotomic character modulo $p$. From this we see that $\mathbf{Q} \subset \mathbf{Q}(\zeta_p) \subset K$. For $p > 2$ this shows that $K$ has no real embeddings as it contains the purely imaginary field $\mathbf{Q}(\zeta_p)$. Let $n = [K : \mathbf{Q}]$ and let $d_K$ be the discriminant of $K/\mathbf{Q}$.

We are now ready for the main theorem of this section.

**Theorem 2.4.** *Let $p \in \{2, 3, 5, 7, 11, 13\}$. Then there are no odd irreducible two dimensional representations of $G = \mathrm{Gal}(\bar{\mathbf{Q}}/\mathbf{Q})$ into $GL_2(\mathbf{F}_p)$ with the following properties: $\rho$ is unramified outside $p$ and is crystalline of Hodge-Tate weight one at $p$.*

We are now ready to prove the theorem.

*Proof.* The cases $p = 2, 3$ were proved by Tate and Serre respectively, without the crystalline hypothesis (see [16] and [14]).

For $p = 5, 7, 11, 13$ we calculate the upper bounds for $|d_K|^{1/n}$ using Theorem 2.1. Next using Diaz y Diaz's tables (see [3]) for Odlyzko bounds on the discriminants get the following upper bounds on $n$–in order of the primes–$n \leq 12, 18, 50, 88$. And in each case the degree must be divisible by $p - 1$.

In each of these cases we see from the above two facts that $K/\mathbf{Q}$ is a tamely ramified extension. Hence by the lemma, we see that

$$|d_K|^{1/n} < p.$$

Hence from Diaz y Diaz's tables we see that $n \leq 4, 6, 24, 40$. In the first two cases, we see that $n = 4, 6$ respectively. For $p = 11$, $n$ can be either 10 or 20 and for $p = 13$, $n$ can be either $12, 24, 36$. Further the extension $K/\mathbf{Q}(\zeta_p)$ must be ramified at the prime lying over $p$–for in each of these case $\mathbf{Q}(\zeta_p)$ has class number one (see [17], page ). One concludes then that in all of these cases $K/\mathbf{Q}$ is totally ramified at $p$ and unramified outside $p$. Moreover the ramification is tame.

As this extension is totally ramified we see that the decomposition group of the (unique) prime lying over $p$ is the whole Galois group. Further, total ramification also implies that there is no residue field extension. Thus the decomposition group is the inertia group. As the extension is tame, all the higher ramification groups are zero. But the inertia group modulo wild inertia subgroup is cyclic (see [15]). Thus the Galois group of $K/\mathbf{Q}$ is Abelian–in fact cyclic. Hence $\rho$ cannot be irreducible and we are done. $\square$

**Remark 2.5.** *One may apply the above result to give a proof of Tate's theorem that there are no elliptic curves over $\mathbf{Z}$.*

*Suppose, if possible, that $E/\mathbf{Q}$ is an elliptic curve with good reduction everywhere. Then by Theorem 2.4. the 5-torsion representation is Abelian and one sees from the proof that all the 25 points of $E[5](\bar{\mathbf{Q}})$ are defined over $\mathbf{Q}(\zeta_5)$. By extension of scalars, $E/\mathbf{Q}(\zeta_5)$ has good reduction at any prime lying over 2. Observe that 2 splits completely in $\mathbf{Q}(\zeta_5)$ and hence the residue field at any prime lying over 2 is $\mathbf{Z}/2$ and hence the elliptic curve has at most $2 + 1 + 2\sqrt{2}$ points. On the other hand, all the 5-torsion injects into the points over $\mathbf{Z}/2$. This is a contradiction.*

*In this case, we observe that all the information we need is contained in the 5-torsion representation–as opposed to general Abelian varieties where we need to prove the results for $p^n$-torsion for suitable $p$ and sufficiently large $n$.*

## 3. Even weights

We can extend the methods slightly and prove that

**Theorem 3.1.** *Let $m \geq 2$ be an integer. Then for $p = 5, 7, 11$ there are no absolutely irreducible representations $\rho$ such that $\rho : G \to GL_m(\mathbf{F}_p)$ is unramified outside $p$ and is crystalline at $p$ and of Hodge-Tate weight at most two.*

*Proof.* It would clearly suffice to show that $\rho \oplus \chi_p$ where $\chi_p$ is the standard cyclotomic character has Abelian image. We observe that the compositum $K' = K\mathbf{Q}(\zeta_p)$ is the fixed field of $\rho \oplus \chi_p$. Moreover $K'$ is a field containing a purely imaginary field $\mathbf{Q}(\zeta_p)$. Thus $K'$ has no real embedding and is also purely imaginary. By the discriminant bound we see that for $p = 5, 7, 11$, we get $n \leq 26, 42, 154$ respectively. And in each case the degree is divisible by $4, 6, 10$. So that $K'/\mathbf{Q}(\zeta_p)$ is an extension of degree $n' \leq 6, 7, 15$ respectively. In each of the case there is one exceptional case when the degree is divisible by the corresponding prime–$n = 20, 42, 110$. Every other possibility gives a tamely ramified extension $K'/\mathbf{Q}$ which is unramified outside $p$.

Thus our tame bound gives $n \leq 6, 10, 24$. In the first two case a simple divisiblity argument gives that $n = 4, 6$ respectively. In the last case we note that $\mathbf{Q}(\zeta_{11})$ has class number one and so $K'/\mathbf{Q}(\zeta_{11})$ which is at most quadratic is also ramified and then $K'/\mathbf{Q}$ is totally and tamely ramified and hence the Galois groups is cyclic and we are done.

We now do the wild case. By [15], the inertia group is a semi direct product of its tame and wild parts, and to handle this case we note, by the lemma below, that the group of order 20 resp. $42, 110$ cannot act irreducibly on a $\mathbf{F}_5$ resp. $\mathbf{F}_7, \mathbf{F}_{11}$ vector space of dimension at least two unless the image its 5-Sylow subgroup (resp. $7, 11$-Sylow subgroup) operates trivially and again we are in the tame case. This completes the proof. □

**Lemma 3.2.** *Let $G$ be a semi-direct product of $(\mathbf{Z}/p)^*$ by $\mathbf{Z}/p$. Suppose $G$ acts on a vector space $V$ of dimension at least two over $\mathbf{F}_p$. Then this representation contains a fixed vector unless image of $\mathbf{Z}/p$ is trivial.*

*Proof.* Let $H$ be the normal subgroup of order $p$. Assume that the image of $H$ under $\rho$ is non-trivial. Then we note that image of any generator of $H$ under $\rho$ is unipotent and if $v$ is an eigenvector for $H$ then $hv = v$ for any $h \in H$. Then as $G$ normalizes $H$ so $hg = gh'$ for some $h' \in H$. And so for any $g \in G$ we have $hgv = gh'v = gv$ so $gv = \lambda v$ for some non-zero $\lambda$. □

## 4. Semi-stable case

In this section we assume the following. Let $\rho : G \to GL_n(\mathbf{F}_p)$ is is a continuous representation which is crystalline at $p$ of Hodge-Tate weights between $[0, p-1]$ and which is semi-stable at primes in a finite set $S$. In this section we use the discriminant bound in the semi-stable case to prove non-existence of irreducible two dimensional representations.

**Theorem 4.1.** *Let $m \geq 2$ be an integer. There are no irreducible, continuous representations $\rho : G \to GL_m(\mathbf{F}_p)$ with either of the following properties:*

1. *Case $S = \{2\}, p = 3$:*
   (a) *$\rho$ is finite flat at 3,*
   (b) *$\rho$ is semi-stable at 2,*
   (c) *$\rho$ is unramified outside $\{2, 3\}$.*
   *or*
2. *Case $S = \{2\}, p = 5$:*
   (a) *$\rho$ is finite flat at 5,*
   (b) *$\rho$ is semi-stable at 2,*
   (c) *$\rho$ is unramified out side $\{2, 5\}$.*

*Proof.* We first prove the case $S = \{2, \}, p = 3$.

We may replace $\rho$ by $\rho \oplus \mu_3$ and assume that $K$ contains $\mathbf{Q}(\zeta_3)$, where as usual $K$ is the fixed field of kernel of $\rho$ and $n = [K : \mathbf{Q}]$. Thus $K$ has non real embeddings. As $K/Q$ is tamely ramified at 2 and possibly wildly ramified at 3. By Theorem 2.1 we see that the discriminant of such an extension satisfies

$$|d_K|^{1/n} < 10.39 \cdots$$

From the tables of Diaz-y-Diaz we see that $n \leq 22$. Note that $K$ is ramified at 2 with ramification index 3 and also ramified at 3 with ramification index divisible by 2. Thus $6|n$ and hence the only possible degrees are $6, 12, 18$. But as 18 does not divide $\#GL_2(\mathbf{F}_3) = 48$ we see that the possible degrees are $6, 12$.

If $n$ is cyclic of degree 6 we are done. If the image is $S_3$ then the representation $\rho \oplus \mu_3$ has no irreducible submodules of dimension bigger than one as $S_3$ has no irreducible representations of dimension bigger than one over $\mathbf{F}_3$. Thus $\rho$ has no irreducible submodules of dimension bigger than one and in particular it is not irreducible.

If $n = 12$, then the only possibilities are that the image is either dihedral group $D_{12}$ or $A_4$. But as $K/\mathbf{Q}(\zeta_3)$, we see that the Galois group of $K/\mathbf{Q}$ has a normal subgroup of order 6; as $A_4$ has no normal subgroup of order 6 we see that the Galois group must be $D_{12}$ in which case the 3-Sylow subgroup of $D_{12}$ is normal. So consider the fixed field $K'$ of the 3-Sylow subgroup. The degree of $K'$ is 4 and as $K/\mathbf{Q}$ is ramified only at $2, 3$ so if $K'$. In fact as the ramification index of any prime lying over 2 in $K$ is 3 we see that 2 is unramified in $K'$. Thus $K'$ is tamely ramified at 3 and unramified outside 3. By Diaz-y-Diaz' tables and the tame discriminant bound we see that $[K' : \mathbf{Q}] \leq 2$. This is a contradiction. Thus the image of $\rho$ cannot be $D_{12}$ and so the only possibility is that $[K : \mathbf{Q}]$ has degree at most 6, and the representation $\rho$ has no irreducible submodules of dimension bigger than 1.

Now we prove the theorem for $S = \{2\}$, $p = 5$.

As before we may assume by replacing $\rho$ by $\rho \oplus \mu_5$ that the fixed field of the kernel of $\rho$ contains $\mathbf{Q}(\zeta_5)$. The discriminant bound gives $n \leq 64$. As $K$ is tamely ramified at 2 with inertia group of order 5 and ramified at 5 with index divisible by 4 we see that $20|n$ and so the only possible degrees for $n$ in this range are $20, 40, 60$. As $K \supset \mathbf{Q}(\zeta_5)$ we see that in each of these cases $\text{Gal}(K/\mathbf{Q})$ has a normal sub group of order $5, 10, 15$ and hence is solvable. As any group of order 15 is cyclic, and as $GL_2(\mathbf{F}_5)$ contains no element of order 15, we see that this group of order 60 is not possible. Thus the Galois group has order 20 or 40.

If the Galois group has order 40, the 5 Sylow subgroup is normal. So passing to its fixed field gives an extension $K'/\mathbf{Q}$ with $[K' : \mathbf{Q}] = 8$ and $K'$ must be unramified at primes lying over 2 because $K/\mathbf{Q}$ is tame at 2. Thus $K'/\mathbf{Q}$ is a tamely ramified extension which is ramified only at 5 and is unramified outside 5. By the tame bound we see that any such extension has degree at most 6 which contradicts the fact that $K'$ has degree 8. Thus the Galois group cannot be of order 40. Thus the only possibility is that the image has order 20. In which case the 5 Sylow subgroup is normal and on every irreducible submodule of $\rho$ this group acts trivially and so any irreducible submodule of $\rho$ is one dimensional. □

## 5. A finiteness result

Let $K$ be a number field. We fix a coefficient field $F_\wp$ which is finite extension of $\mathbf{Q}_p$. We will say a continuous, irreducible representation $\rho : \text{Gal}(\bar{K}/K) \to GL_n(F_\wp)$ is "effective of motivic type" if the following conditions are satisfied.

1. $\rho$ is unramified outside a finite set of primes of $K$,
2. $\rho$ has weight $k \geq 1$ in the sense that for all primes $\mathfrak{q}$ for which $\rho$ is unramified, the eigenvalues of $\text{Frob}_\mathfrak{q}$ are algebraic integers of absolute value $N(\mathfrak{q})^{k/2}$, where $N(\mathfrak{q})$ is the norm of the ideal $\mathfrak{q}$.

REMARKS ON METHODS OF FONTAINE AND FALTINGS 7It is standard (see [2]) that representations arising from algebraic geometry are of this type. The method of proof of the following result is extracted from Faltings' proof of the Shafarevitch conjecture (see [4] and [5]).

**Theorem 5.1.** *Fix integers $n \geq 1, k \geq 1$, and a finite set of primes $S$. Then the set of isomorphism classes of continuous, irreducible representations $\rho : \mathrm{Gal}(\bar{K}/K) \to GL_n(F_\wp)$ which are unramified outside $S$ and which are effective of motivic type of weight $k$, is finite.*

*Proof.* The proof given in Chapter V, Section 2 of [5] adapts with out any difficulty to the above situation. We will briefly recall his proof here for the reader's convenience. We need a version of Chebotarev density theorem (see Chapter V, Corollary 2.4; we also need a the Hermite-Minkowski Theorem (see Theorem 2.6 of Chapter V).

Next fix a prime $\mathfrak{q} \notin S$ of $K$ not lying over any prime lying over $p$ in $K$; also fix a representation $\rho$ which is effective of motivic type which is unramified outside $S$. Now we note that there are only finitely many possibilities for the local $L$-factor $\det(1 - N(\mathfrak{q})^{-s}\rho(\mathrm{frob}_\mathfrak{q}))$ of $\rho$ at $\mathfrak{q}$. This follows (see Chapter V, Lemma 2.6) from the fact that there are only finitely many algebraic numbers with bounded conjugates and bounded degree.

Let $\tilde{K}$ be a finite Galois extension which contains all extensions of $K$ which are finite and Galois over $K$ and which are unramified outside $S$ and which have degree at most $p^{2[F_\wp:\mathbf{Q}_p]n^2}$ over $K$. Let $T$ be the set of primes whose existence is guaranteed by Chebotarev density theorem (Chapter V, Corollary 2.4) for the extension $\tilde{K}/K$.

The main lemma of Faltings' (see Proposition 2.7 of Chapter V) is that if $\rho_1, \rho_2$ are two irreducible representations of the above type such that the traces of Frobenius coincide on the set of primes $T$ constructed above then, $\rho_1 \simeq \rho_2$. The proof of this Proposition as given in [5] works in our case too. Now the proof of the above theorem can be completed as the proof of Theorem 2.8 of Chapter V. □

It follows from the above theorem that the set of irreducible representations which arise from algebraic geometry of fixed ramification data, weight, dimension and coefficient field is finite. One of curious consequences of this is the following;

**Corollary 5.2.** *Fix a finite set $S$ of primes of $\mathbf{Q}$. Fix an extension $F_\wp$ of $\mathbf{Q}_p$. Fix an integer $k \geq 2$. Then there are only finite number of new forms of weight $k$ and of level divisible only by primes in $S$ and such that corresponding two dimensional representations have values in $GL_2(F_\wp)$. In particular, there are only finitely many irreducible two dimensional representations $\rho : \mathrm{Gal}(\bar{\mathbf{Q}}/\mathbf{Q}) \to GL_2(\mathbf{Q}_p)$ which are unramified outside $p$ and which are of weight 11 and such that the eigenvalues of Frobenius elements are algebraic integers.*

## References

[1] Sharon Brueggeman, *The nonexistence of certain Galois representations unramified outside 5*. To appear in Journal of Number Theory.

[2] P. Deligne, *La conjecture de Weil* I Inst. Hautes Études Sci. Publ. Math., **43**, (1974), pp 273-307.

[3] F. Diaz y Diaz, *Tables minorant la racine n-ième du discriminant d'un corps de degré n*. Ph. D Thesis, Publ. Math. d'Université d'Orsay, 1980.

[4] Gerd Faltings, *Endlichkeitssätze für abelsche Variétés über zahlkörpern.* Invent. Math. **73**, (1983), pp 349-366.

[5] Gerd Faltings, Gisbert Wūstholz et al. *Rational points*, Aspects of mathematics, E6, F. Vieweg & Sohn.
[6] J.-M. Fontaine, *Sur certains types de représentations p-adiques du groupe de Galois d'un corps local; construction d'un anneau de Barsotti-Tate*, Ann. of Math. **115** (1982), no. 3, 529–577.
[7] ______, *Il n'y a pas de variété abélienne sur* **Z**, Invent. Math. **81** (1985), 515–538.
[8] ______, *Schémas propres et lisses sur* **Z**, Proceedings of the Indo-French Conference on Geometry (Delhi), Hindustan Book Agency, 1993, pp. 43–56.
[9] J.-M. Fontaine and L. Laffaille, *Construction de représentations p-adiques.*, Ann. Sci. cole Norm. Sup. **15** (1983), no. 4, 547–608.
[10] J.-M. Fontaine and B. Mazur, *Geometric Galois representations I*, Elliptic curves, modular forms and Fermat's last theorem (Hong Kong 1993) (Cambridge, MA), Number Theory, no. 1, International Press, pp. 41–78.
[11] John. W. Jones and David .P. Roberts, *Sextic number fields with discriminant* $-^j 2^a 3^b$, Centre de Recherches de Mathematiques, CRM proceedings and Lecture Notes in Mathematics, Vol 18, 1998.
[12] G. Poitou, *Minorations de discriminants (d'après a. m. Odlyzko)*, Séminaire Bourbaki, vol. 1975-76, Lecture Notes in Mathematics, no. 567, Springer, pp. 136–153.
[13] Kenneth Ribet, *On modular representations of* $\text{Gal}(\overline{\mathbf{Q}}/\mathbf{Q})$ *arising from modular forms*, Invent. Math. **100** (1990), no. 2, 431–476.
[14] Jean-Pierre Serre, *Sur les représentations modulaires de degré 2 de* $\text{Gal}(\overline{\mathbf{Q}}/\mathbf{Q})$., Duke Math. J. **54** (1987), no. 1, 179–230.
[15] ______, *Local fields*, Graduate Texts in Mathematics, vol. 67, Springer-Verlag, New York-Berlin, 1979.
[16] John Tate, *The non-existence of certain Galois extensions of* **Q** *unramified outside* 2, Arithmetic geometry (Providence, RI), Contemp. Math., vol. 174, Amer. Math. Soc., 1994, pp. 153–156.
[17] L. Washington, *Introduction to cyclotomic fields*, second edition ed., Graduate Texts in Mathematics, vol. 83, Springer-Verlag, 1997.

Department of Mathematics, University of Arizona, 617 N Santa Rita, P O Box 210089, Tucson, AZ 85721, USA

*E-mail address*: `kirti@math.arizona.edu`